\documentclass{svjour3}                     
\smartqed  
\usepackage{graphicx}
\usepackage{amssymb,amsmath,amsfonts,latexsym,euscript}
\begin{document}

\title{The order of convergence of an optimal quadrature formula with derivative in the space $W_2^{(2,1)}$}


\titlerunning{The order of convergence of an optimal quadrature formula}        

\author{Abdullo R. Hayotov$^*$, Rashidjon G. Rasulov}


\institute{A.R. Hayotov \at
$^*$Corresponding author:
Department of Mathematical Sciences, KAIST, 291 Daehak-ro, Yuseong-gu, Daejeon 34141, Republic of Korea\\
\email{hayotov@kaist.ac.kr}.\\
V.I.Romanovskiy Institute of Mathematics, Uzbekistan Academy of Sciences,
81 M.Ulugbek str., Tashkent 100170, Uzbekistan,\\
\email{hayotov@mail.ru}.
\and
 R.G. Rasulov \at
V.I.Romanovskiy Institute of Mathematics, Uzbekistan Academy of Sciences,
81 M.Ulugbek str., Tashkent 100170, Uzbekistan\\
\email{r.rasulov1990@mail.ru}.
}

\date{Received: date / Accepted: date}

\maketitle
\begin{abstract}
The present work is devoted to extension of the trapezoidal rule in the space $W_2^{(2,1)}$.
The optimal quadrature formula is obtained by minimizing the error of the formula by coefficients
at values of the first derivative of a integrand.
Using the discrete analog of the operator $\frac{d^2}{dx^{2}}-1$ the explicit formulas for the coefficients of the
optimal quadrature formula are obtained. Furthermore, it is proved that the obtained quadrature formula is exact for any function of the set $\mathbf{F}=\mathrm{span}\{1,x,e^{x},e^{-x}\}$. Finally, in the space $W_2^{(2,1)}$ the square of the norm of the error
functional of the constructed quadrature formula is calculated. It is shown that the error of the obtained optimal quadrature formula
is less than the error of the Euler-Maclaurin quadrature formula on the space $L_2^{(2)}$.

\textbf{MSC:} 65D30, 65D32.

\textbf{Keywords:} optimal quadrature formula, Hilbert space, the error functional, S.L.~Sobolev's method, discrete
argument function, the order of convergence.
\end{abstract}

\maketitle

\section{Introduction}

It is known, that quadrature and cubature formulas, are methods for the approximate evaluation of definite integrals.
In addition and even more important, quadrature formulas provide a basic and important tool for the
numerical solution of differential and integral equations. The theory of cubature formulas consists mainly of three branches
dealing with exact formulas, formulas based on functional-analytic methods, and formulas based on probabilistic methods \cite{Sobolev74,SobVas}.
In the functional-analytic methods the error between
an integral and corresponding cubature sum is considered as a linear functional on a Banach space and it is estimated
by the norm of the error functional in the conjugate Banach space.
The norm of the error functional depends on coefficients
and nodes of the formula. The problem of finding the minimum of the norm of the error functional by coefficients and by nodes is called \emph{S.M.Nikol'skii problem}, and the obtained formula is called \emph{the optimal formula in the sense of Nikol'skii} (see, for instance, \cite{Nik88}). Minimization of the norm of the error functional by coefficients when the nodes are fixed is called \emph{Sard's problem}.
And the obtained formula is called \emph{the optimal formula in the sense of Sard}. First this problem was studied by A. Sard \cite{Sard}.
Solving these problems in different spaces of differentiable functions various type of optimal formulas of numerical integration are obtained.

There are several methods for constructing the optimal quadrature
formulas in the sense of Sard such as  the spline method, the
$\phi-$ function method (see e.g. \cite{BlaCom},
\cite{SchSil}) and the Sobolev method. It should be noted that the
Sobolev method is based on using a discrete analog
of a linear differential operator (see e.g. \cite{Sobolev06,Sobolev74,SobVas}). In different spaces based on
these methods, the Sard problem was studied by many authors,
see, for example, \cite{IBab,BlaCom,CatCom,HayMilShad10,HayMilShad15,Koh,FLan,Sard,SchSil,ShadHay11,Sobolev06,Sobolev74,SobVas,Zag}
and references therein.

Among these formulas the Euler-Maclaurin type quadrature formulas are very important for numerical integration of differentiable functions
and are widely used in applications. In different spaces the optimality of the Euler-Maclaurin type
quadrature and cubature formulas were studied, for instance, in works \cite{CatCom,FLan,Mic74,Schoen65,ShadHayNur13,ShadHayNur16,ShadNur18,Zhen81}.

The Euler-Maclaurin quadrature formulas can be viewed as well as an extension of the trapezoidal rule by the inclusion of correction terms.
It should be noted that in applications and in solution of practical problems numerical integration formulas are interesting for functions with small smoothness.

The present paper is also devoted to extension of the trapezoidal rule.

We consider a quadrature formula of the form
\begin{equation}\label{(1.1)}
\int\limits_0^1 \varphi (x)dx \cong  \sum\limits_{\beta  = 0}^N
(C_0[\beta]\varphi(h\beta)+C_1[\beta]\varphi'(h\beta))
\end{equation}
where  $C_0[\beta]$ are coefficients of the trapezoidal rule, i.e.
\begin{equation}\label{(C0)}
\begin{array}{lll}
C_0[0]&=&\frac{h}{2},\\
C_0[\beta]&=&h,\ \beta=1,2,...,N-1,\\
C_0[N]&=&\frac{h}{2},
\end{array}
\end{equation}
$C_1[\beta]$ are unknown coefficients of the formula (\ref{(1.1)}) and they should be found, $h=\frac{1}{N}$, $N$ is a natural number.
We suppose that an integrand $\varphi$ belongs to $W_2^{(2,1)}(0,1)$,
where by $W_2^{(2,1)}(0,1)$ we denote the class of all functions $\varphi$ defined on $[0,1]$ which posses an absolutely continuous
first derivative and whose second derivative is in $L_2(0,1)$. The class $W_2^{(2,1)}(0,1)$ under the pseudo-inner product
$$
\langle\varphi,\psi\rangle=\int\limits_0^1(\varphi''(x)+\varphi'(x))(\psi''(x)+\psi'(x))dx
$$
is a Hilbert space if we identify functions that differ by a linear combination of a constant and $e^{-x}$ (see, for example, \cite{Ahlb67}).
Here, in the Hilbert space $W_2^{(2,1)}(0,1)$, we consider the corresponding norm
\begin{equation}\label{(1.3)}
\|\varphi|W_2^{(2,1)}(0,1)\|=\left[\int\nolimits_0^1(\varphi''(x)+\varphi'(x))^2dx\right]^{1/2}.
\end{equation}

The difference
\begin{equation}\label{(1.6)}
(\ell,\varphi) = \int\limits_0^1 \varphi (x)dx -  \sum\limits_{\beta  = 0}^N
(C_0[\beta]\varphi(h\beta)+C_1[\beta]\varphi'(h\beta))
\end{equation}
is called \emph{the error} and
\begin{equation}\label{(1.2)}
\ell (x) = \varepsilon _{[0,1]} (x) - \sum\limits_{\beta  = 0}^N
(C_0[\beta]\delta(x - h\beta)-C_1[\beta]\delta'(x-h\beta)),
\end{equation}
is said to be \emph{the error functional} of the quadrature formula (\ref{(1.1)}), where $\varepsilon
_{[0,1]} (x)$ is the indicator of the interval $[0,1]$ and $\delta$ is Dirac's delta function.
The value $(\ell,\varphi)$ of the error functional $\ell$ at a function $\varphi$ is defined as
$(\ell,\varphi)=\int_{-\infty}^{\infty}\ell(x)\varphi(x)dx.$

In order that the error functional (\ref{(1.2)}) is defined on the space $W_2^{(2,1)}(0,1)$ it is necessary to impose the following conditions
for the functional $\ell$
\begin{eqnarray}
(\ell,1):=&&1-\sum\limits_{\beta=0}^NC_0[\beta]=0,\label{(1.4)}\\
(\ell,e^{-x}):=&&\int\limits_0^1e^{-x}dx-\sum\limits_{\beta=0}^N(C_0[\beta]e^{-h\beta}-C_1[\beta]e^{-h\beta})=0. \label{(1.5)}
\end{eqnarray}
The last two equations mean that the quadrature formula (\ref{(1.1)}) is exact for any constant and $e^{-x}$.
We have chosen the coefficients $C_0[\beta]$, $\beta=0,1,...,N$ such that the equality (\ref{(1.4)}) is fulfilled. Therefore we have only condition (\ref{(1.5)}) for
coefficients $C_1[\beta]$, $\beta=0,1,...,N$.

The error functional $\ell$ of the formula (\ref{(1.1)}) is a linear functional
in $W_2^{(2,1)*}(0,1)$, where $W_2^{(2,1)*}(0,1)$ is the
conjugate space to the space $W_2^{(2,1)}(0,1)$.

By the Cauchy-Schwarz inequality we have the following
$$
|(\ell,\varphi)|\leq \|\varphi|W_2^{(2,1)}(0,1)\|\cdot
\|\ell|W_2^{(2,1)*}(0,1)\|.
$$
Hence we conclude that the error (\ref{(1.6)}) of the formula (\ref{(1.1)}) is
estimated by the norm
\begin{equation}\label{norm}
\left\| {\ell |W_2^{(2,1)*}(0,1) } \right\| = \mathop {\sup
}\limits_{\left\| {\varphi |W_2^{(2,1)}(0,1)} \right\| = 1}
\left| {\left( {\ell ,\varphi } \right)} \right|
\end{equation}
of the error functional (\ref{(1.2)}).

The main aim of this work is to find the minimum of the absolute value of the error (\ref{(1.6)}) by coefficients
$C_1[\beta]$ for given $C_0[\beta]$ in the space $W_2^{(2,1)}$.
That is the problem is to find the coefficients $C_1[\beta]$ that satisfy the following
equality
\begin{equation}\label{(1.7)}
\left\| \mathring{\ell} |W_2^{(2,1)*}
\right\| = \mathop {\inf }\limits_{C_1[\beta] } \left\| \ell |W_2^{(2,1)*} \right\|.
\end{equation}
The coefficients $C_1[\beta]$ which satisfy the last equation are called \emph{optimal} and are denoted as $\mathring{C}_1[\beta]$.

Thus, to obtain the optimal quadrature formula of the form (\ref{(1.1)}) in the sense of
Sard in the space $W_2^{(2,1)}(0,1)$, we need to solve the
following problems.

\medskip

\textbf{Problem 1.} \emph{Find the norm of the error functional
 (\ref{(1.2)}) of the quadrature formula (\ref{(1.1)}) in the space $W_2^{(2,1)*}$.}

\medskip

\textbf{Problem 2.} \emph{Find the coefficients $\mathring{C}_1[\beta] $ that
satisfy equality (\ref{(1.7)}).}

\medskip

Here we solve Problems 1 and 2 by Sobolev's method using the discrete analog
of the differential operator $\frac{d^2}{dx^{2}}-1$.

The paper is organized as follows: in Section 2 using
the extremal function of the error functional $\ell$  the
norm of this functional is calculated, i.e. Problem 1 is solved; Section 3 is devoted to solution of Problem 2. Here  the system
of linear equations for the coefficients $C_1[\beta]$ of the optimal quadrature
formulas (\ref{(1.1)}) is obtained in the space $W_2^{(2,1)}(0,1)$. In Subsection 3.1 using the discrete analog of the operator $\frac{d^2}{dx^{2}}-1$ the explicit formulas for the coefficients $C_1[\beta]$ of optimal quadrature formula of the form (\ref{(1.1)}) are obtained.
Furthermore, it is proved that the obtained quadrature formula of the form (\ref{(1.1)}) is exact for any function of the set $\mathbf{F}=\mathrm{span}\{1,x,e^{x},e^{-x}\}$.
Finally, in Subsection 3.2 in the space $W_2^{(2,1)}$ the square of the norm of the error
functional of the constructed quadrature formula is calculated. It is shown that the error of the obtained optimal quadrature formula
is less than the error of the Euler-Maclaurin quadrature formula on the space $L_2^{(2)}$.

\section{The norm of the error functional (\ref{(1.2)})}

In this section we study Problem 1. To calculate the norm of the error
functional (\ref{(1.2)}) in the space $W_2^{(2,1)*}(0,1)$
we use \emph{the extremal function} for this functional (see, \cite{Sobolev74,SobVas})  which satisfies the equality
\begin{equation*}
 \left( \ell,\psi _\ell \right) = \left\| {\ell
|W_2^{(2,1)*}} \right\| \cdot \left\| {\psi
_\ell|W_2^{(2,1)}} \right\|.
\end{equation*}

We note that in \cite{ShadHay14} for a linear functional $\ell$ defined on the Hilbert space $W_2^{(m,m-1)}$ the extremal function was found
and it was shown that the extremal function $\psi_{\ell}$ is the solution of the boundary value problem
\begin{eqnarray}
&& \psi_{\ell}^{(2m)}(x)-\psi_{\ell}^{(2m-2)}(x)=(-1)^m\ell(x),\label{(2.1)}\\
&&\left(\psi_{\ell}^{(m+s)}(x)-\psi_{\ell}^{(m+s-2)}(x) \right)|_{x=0}^{x=1}=0,\ \ s=0,1,...,m-1,\label{(2.2)}\\
&&\left(\psi_{\ell}^{(m)}(x)+\psi_{\ell}^{(m-1)}(x) \right)|_{x=0}^{x=1}=0.\label{(2.3)}
\end{eqnarray}
That is for the extremal function $\psi_{\ell}$ the following was proved.

\begin{theorem}[Theorem 2.1 of \cite{ShadHay14}]\label{THM2.1} The solution of the boundary value problem (\ref{(2.1)})-(\ref{(2.3)}) is the extremal function
$\psi_{\ell}$ of the error functional $\ell$ and has the following form
$$
\psi_{\ell}(x)=(-1)^m\ell(x)*G_m(x)+P_{m-2}(x)+de^{-x},
$$
where
\begin{equation}\label{Gm}
G_m(x)=\frac{\mathrm{sgn}x}{2}\left(\frac{e^x-e^{-x}}{2}-\sum\limits_{k=1}^{m-1}\frac{x^{2k-1}}{(2k-1)!}\right)
\end{equation}
is the solution of the eqution
$
G_m^{(2m)}(x)-G_m^{(2m-2)}(x)=\delta(x),
$
$d$ is any real number and $P_{m-2}(x)$ is a polynomial of degree $m-2$.
\end{theorem}

Furthermore, there were shown that
$\|\ell|W_2^{(m,m-1)}\|=\|\psi_{\ell}|W_2^{(m,m-1)*}\|$
and
\begin{equation}\label{(2.4)}
\left( \ell,\psi _\ell \right) = \left\| \ell
|W_2^{(m,m-1)*} \right\|^2.
\end{equation}

From Theorem \ref{THM2.1}, in the case $m=2$,   we get the extremal function $\psi_{\ell}$ for the error functional (\ref{(1.2)}) and it has the form
\begin{equation}\label{(2.5)}
\psi_{\ell}(x)=\ell(x)*G_2(x)+p_0+de^{-x},
\end{equation}
where
\begin{equation}\label{(2.6)}
G_2(x)=\frac{\mathrm{sgn}x}{2}\left(\frac{e^x-e^{-x}}{2}-x\right),
\end{equation}
$d$ and $p_0$ are any real numbers.

Then, from (\ref{(2.4)}), in the case $m=2$, using (\ref{(1.2)}) and (\ref{(2.5)}), taking into account equations
(\ref{(1.4)}) and (\ref{(1.5)}),
we get
\begin{equation}\label{(2.7)}
\begin{array}{l}
\|\ell\|^2=(\ell,\psi_{\ell})=\sum\limits_{\beta=0}^N\sum\limits_{\gamma=0}^N\bigg(C_0[\beta]C_0[\gamma]G_2(h\beta-h\gamma)-C_1[\beta]C_1[\gamma]G''_2(h\beta-h\gamma)\bigg)+\\
\qquad\qquad \qquad\qquad +2\sum\limits_{\beta=0}^NC_1[\beta]\bigg(\int\limits_0^1G'_2(x-h\beta)dx+\sum\limits_{\gamma=0}^NC_0[\gamma]G'_2(h\beta-h\gamma)\bigg)-\\
\qquad\qquad \qquad\qquad -2\sum\limits_{\beta=0}^NC_0[\beta]\int\limits_0^1G_2(x-h\beta)dx+\int\limits_0^1\int\limits_0^1G_2(x-y)dxdy,\\
\end{array}
\end{equation}
where $G_2(x)$ is defined by (\ref{(2.6)}), $G_2'(x)$ and $G_2''(x)$ are derivatives of $G_2(x)$, i.e.
\begin{equation}\label{(2.8)}
G'_2(x)=\frac{\mathrm{sgn} x}{2}\left(\frac{e^x+e^{-x}}{2}-1\right)\mbox{ and } G''_2(x)=\frac{\mathrm{sgn} x}{2}\left(\frac{e^x-e^{-x}}{2}\right).
\end{equation}
It is easy to see from (\ref{Gm}) and (\ref{(2.8)}) that
\begin{equation}\label{(2.10)}
G_1(x)=G_2''(x).
\end{equation}
Thus Problem 1 is solved.

In the next section we study Problem 2.

\section{Minimization of the norm (\ref{(2.7)})}

Now we consider the minimization problem  of the expression (\ref{(2.7)}) by the coefficients $C_1[\beta]$ under the condition
(\ref{(1.5)}). For this we use the Lagrange method of conditional extremum.

Consider the Lagrange function
\begin{equation*}
\Psi(C_1[0],C_1[1],...,C_1[N],d)=\|\ell\|^2+2d(\ell,e^{-x}).
\end{equation*}
Taking partial derivatives from the function $\Psi$  by $C_1[\beta]$, $\beta=0,1,...,N$  then equating them to 0
and using the condition (\ref{(1.5)}), we get the following system of $N+2$ linear equations with $N+2$ unknowns
\begin{eqnarray}
&&\sum\limits_{\gamma=0}^NC_1[\gamma]G''_2(h\beta-h\gamma)+de^{-h\beta}=F(h\beta),\ \ \beta=0,1,2,...,N,\label{(3.1)}\\
&& \sum\limits_{\gamma=0}^NC_1[\gamma]e^{-h\gamma}=g,\label{(3.2)}
\end{eqnarray}
where
\begin{eqnarray}
F(h\beta)&=&\int_0^1G'_2(x-h\beta)dx+\sum\limits_{\gamma=0}^NC_0[\gamma]G'_2(h\beta-h\gamma),\label{(3.3)}\\
g&=&e^{-1}-1+\sum\limits_{\gamma=0}^NC_0[\gamma]e^{-h\gamma}.\label{(3.4)}
\end{eqnarray}
Here $C_0[\gamma]$, $\gamma=0,1,...,N$ are defined by (\ref{(C0)}), $C_1[\beta]$, $\beta=0,1,...,N$ and $d$ are unknowns.

The system (\ref{(3.1)})-(\ref{(3.2)}) has a unique solution for any fixed natural number $N$ and this solution
gives the minimum to the expression (\ref{(2.7)}). Here we
omit the proof of the existence and uniqueness of the solution of this system. These statements can be proved
similarly as the proof of the existence and uniqueness of the solution of the discrete Wiener-Hopf type system for the optimal coefficients  of quadrature formulas with the form $\int_0^1f(x)dx\cong \sum_{\beta=0}^NC[\beta]f[\beta]$ in
the space $L_2^{(m)}(0,1)$ (see \cite{Sobolev06,Sobolev74,SobVas}).

\medskip

\textbf{3.1. The coefficients of the optimal quadrature formula (\ref{(1.1)})}

\medskip

In this subsection we solve the system (\ref{(3.1)})-(\ref{(3.2)}) and find the explicit forms for optimal coefficients
$C_1[\beta]$, $\beta=0,1,...,N$.

Here we use the concept of discrete argument functions and operations on them. The theory of
discrete argument functions is given in \cite{Sobolev74,SobVas}. We give some definitions about
functions of discrete argument.

Suppose that $\varphi$ and $\psi$ are real-valued
functions of real variable $x$ and are defined in real line
$\mathbb{R}$.

\medskip

A function $\varphi (h\beta)$ is called
\emph{a function of discrete argument} if it is defined on some set
of integer values of $\beta$.

\medskip

\emph{The inner product} of two
discrete functions $\varphi(h\beta )$ and $\psi (h\beta )$ is defined as the
following number
$$
\left[ {\varphi(h\beta),\psi(h\beta) } \right] =
\sum\limits_{\beta  =  - \infty }^\infty  {\varphi (h\beta ) \cdot
\psi (h\beta )},
$$
if the series on the right hand side of the last equality
converges absolutely.

\medskip

\textit{The convolution} of two
discrete argument functions $\varphi(h\beta )$ and $\psi (h\beta )$ is the
inner product
$$
\varphi (h\beta )*\psi (h\beta ) = \left[ {\varphi (h\gamma ),\psi
(h\beta  - h\gamma )} \right] = \sum\limits_{\gamma  =
- \infty}^\infty  {\varphi (h\gamma ) \cdot \psi (h\beta  - h\gamma )}.
$$

\medskip

Furthermore, for finding the coefficients $C_1[\beta]$ of the optimal quadrature formula (\ref{(1.1)})
we need the discrete analog of the differential operator $\frac{d^2}{dx^{2}}-1$. It should be noted that in the work \cite{ShadHay04}
 the discrete analog of the differential operator $\frac{d^{2m}}{dx^{2m}}-\frac{d^{2m-2}}{dx^{2m-2}}$ was constructed.
In particular, when $m=1$ from the result of the work \cite{ShadHay04} we get the following

\begin{theorem}\label{Thm3.1}
The discrete analog $D_1(h\beta)$ of the differential operator $\frac{d^2}{dx^{2}}-1$ satisfying
the equation
$$
D_1(h\beta)*G_1(h\beta)=\delta_d(h\beta)
$$
has the form
\begin{equation}\label{(3.5)}
D_1(h\beta)=
\frac{1}{1-e^{2h}}
\left\{
\begin{array}{ll}
0,& |\beta|\geq 2,\\[2mm]
-2e^h,& |\beta|=1,\\[2mm]
2(1+e^{2h}),& \beta=0,
\end{array}
\right.
\end{equation}
where
$G_1(h\beta)=\frac{\mathrm{sgn}(h\beta)}{2}\left(\frac{e^{h\beta}-e^{-h\beta}}{2}\right)$ and
$\delta_d(h\beta)=
\left\{
\begin{array}{ll}
1,& \beta=0,\\
0,&\beta\neq 0.
\end{array}
\right.
$
\end{theorem}

Furthermore, it is easy to check that
\begin{equation}\label{(3.6)}
D_1(h\beta)*e^{h\beta}=0 \mbox{ and } D_1(h\beta)*e^{-h\beta}=0.
\end{equation}

Now we turn to get the solution of the system (\ref{(3.1)})--(\ref{(3.2)}) using (\ref{(3.5)}).

Suppose $C_1[\beta]=0$ when $\beta<0$ and $\beta>N$.
Then we rewrite the system (\ref{(3.1)})--(\ref{(3.2)}) in the following convolution form
\begin{eqnarray}
&&C_1[\beta]*G''_2(h\beta)+de^{-h\beta}=F(h\beta),\ \ \beta=0,1,...,N,\label{(3.7)}\\
&& \sum\limits_{\gamma=0}^NC_1[\gamma]e^{-h\gamma}=g.\label{(3.8)}
\end{eqnarray}
Here, calculating the right hand sides of equalities (\ref{(3.3)}) and (\ref{(3.4)}), for $F(h\beta)$ and $g$ we get
\begin{eqnarray}
F(h\beta)&=&\left(\frac{e^{h\beta}}{8}(e^{-1}+1)-\frac{e^{-h\beta}}{8}(e+1)\right)\left(\frac{h(e^h+1)}{e^h-1}-2\right),\label{(3.9)}\\
g&=&\frac{1}{2}(1-e^{-1})\left(\frac{h(e^h+1)}{e^h-1}-2\right).\label{(3.10)}
\end{eqnarray}

Taking into account (\ref{(2.8)}) and (\ref{(2.10)}), using Theorem \ref{Thm3.1} we get
\begin{equation}\label{(3.11)}
D_1(h\beta)*G_2''(h\beta)=\delta_d(h\beta),
\end{equation}

Denoting by
\begin{equation}\label{(3.12)}
u(h\beta)=C_1[\beta]*G_2''(h\beta)+de^{-h\beta}
\end{equation}
the left hand side of the equation (\ref{(3.7)}) we get
\begin{equation}\label{(3.13)}
C_1[\beta]=D_1(h\beta)*u(h\beta).
\end{equation}
Indeed, if the function $u(h\beta)$ is defined at all integer values of $\beta$, then taking into account Theorem \ref{Thm3.1} and using properties
(\ref{(3.6)}) of the function $D_1(h\beta)$, we have
\begin{eqnarray*}
D_1(h\beta)*u(h\beta)&=&D_1(h\beta)*(G_2''(h\beta)*C_1[\beta])+D_1(h\beta)*(d\ e^{-h\beta}) \\
                     &=&C_1[\beta]*(D_1(h\beta)*G_2''(h\beta))\\
                     &=&C_1[\beta]*\delta_d(h\beta)\\
                     &=&C_1[\beta].\\
\end{eqnarray*}
Thus, if we find the function $u(h\beta)$ for all integer values of $\beta$ then the optimal coefficients $C_1[\beta]$ will be found
from the equality (\ref{(3.13)}).

The following is true.

\begin{theorem}\label{Thm3.2}
The coefficients of the optimal quadrature formula of the form (\ref{(1.1)}) in the sense of Sard in the space $W_2^{(2,1)}(0,1)$
have the following form
\begin{eqnarray}
\mathring{C}_1[0]&=&\frac{h(e^{h}+1)}{2(e^{h}-1)}-1,\nonumber\\
\mathring{C}_1[\beta]&=&0,\ \ \    \beta=1,2,...,N-1,\label{(3.14)}\\
\mathring{C}_1[N]&=&1-\frac{h(e^{h}+1)}{2(e^{h}-1)}.\nonumber
\end{eqnarray}
\end{theorem}

\textit{Proof.} From equality (\ref{(3.7)}) taking into account (\ref{(3.12)}) we get that
$$
u(h\beta)=F(h\beta)\mbox{ for }\beta=0,1,...,N.
$$
Now we find the function $u(h\beta)$ for $\beta<0$ and $\beta>N$.

Let $\beta\leq 0$ then from (\ref{(3.12)}), using the form (\ref{(2.8)}) of the function $G_2''(x)$ and equality (\ref{(3.8)}), we have
$$
u(h\beta)=-\frac{1}{4}e^{h\beta} g+e^{-h\beta}\frac{1}{4}\sum_{\gamma=0}^{N}C_{1}[\gamma]e^{h\gamma}+de^{-h\beta}.
$$
Similarly, for $\beta\geq N$ we obtain
$$
u(h\beta)=\frac{1}{4}e^{h\beta} g-e^{-h\beta}\frac{1}{4}\sum_{\gamma=0}^{N}C_{1}[\gamma]e^{h\gamma}+de^{-h\beta}.
$$
Then, keeping in mind the last two equalities and denoting by
$$
D=\frac{1}{4}\sum\limits_{\gamma=0}^NC_1[\gamma]e^{h\gamma},
$$
for $u(h\beta)$ we get the following
\begin{equation}\label{(3.15)}
u(h\beta)=\left\{
\begin{array}{ll}
-\frac{1}{4}e^{h\beta} g+(d+D)e^{-h\beta},&\beta\leq 0,\\[2mm]
F(h\beta),&0\leq \beta\leq N,\\[2mm]
\frac{1}{4}e^{h\beta} g+(d-D)e^{-h\beta}, &\beta\geq N.
\end{array}
\right.
\end{equation}
Here, in the equality (\ref{(3.15)}), $d$ and $D$  are unknowns.
These unknowns can be found from the conditions of consistency of values of the function $u(h\beta)$
at the points $\beta=0$ and $\beta=N$. Therefore from (\ref{(3.15)}) when $\beta=0$ and $\beta=N$ we obtain the system of linear equations for
$d$ and $D$. Then, using (\ref{(3.9)}) and (\ref{(3.10)}), after some calculations, we have
\begin{equation}\label{(3.16)}
\left\{
\begin{array}{l}
d=0,\\[2mm]
D=\frac{1}{8}\left(\frac{h(e^h+1)}{e^h-1}-2\right)(1-e).
\end{array}
\right.
\end{equation}
Finally, from (\ref{(3.13)}) for $\beta=0,1,...,N$, using (\ref{(3.5)}) and (\ref{(3.15)}) and taking into account (\ref{(3.16)}),
by directly calculations
we get (\ref{(3.14)}). Theorem \ref{Thm3.2} is proved. \hfill $\Box$

\medskip

\textbf{Remark 1.} Using (\ref{(C0)}) and (\ref{(3.14)}),
one can get that $(\ell,x)=0$ and $(\ell,e^x)=0$. These equalities mean that the optimal quadrature formula of the form (\ref{(1.1)})
with the coefficients (\ref{(C0)}) and (\ref{(3.14)}) is also exact for functions $x$ and $e^x$. Therefore, keeping in mind equalities (\ref{(1.4)}) and (\ref{(1.5)}), we conclude that the optimal quadrature formula of the form (1.1) with  coefficients (\ref{(C0)}) and (\ref{(3.14)}) is exact for any linear combinations of functions $1,\ x,\ e^x$ and $e^{-x}$, i.e. it is exact for elements of the set $\mathbf{F}=\mathrm{span}\{1,\ x,\ e^x,\ e^{-x}\}$.

\medskip

\textbf{3.2. The norm of the error functional $\mathring\ell$ of the optimal quadrature formula (\ref{(1.1)})}

\medskip

In this subsection we study the order of convergence of the optimal quadrature formula of the form (\ref{(1.1)})
with coefficients (\ref{(C0)}) and (\ref{(3.14)}), i.e. we calculate the square of the norm (\ref{(2.7)}) of the error functional for the optimal quadrature formula
(\ref{(1.1)}).

The following holds

\begin{theorem}\label{Thm3.3}
Square of the norm of the error functional (\ref{(1.2)}) for the optimal quadrature formula (\ref{(1.1)}) with coefficients
(\ref{(C0)}) and (\ref{(3.14)}) on the space $W_2^{(2,1)}(0,1)$ has the form
\begin{eqnarray}
\left\|\mathring{\ell}|W_2^{(2,1)*}(0,1)\right\|^2&=&-\sum\limits_{n=4}^{\infty}\frac{B_nh^n}{n!}\nonumber\\
&=&\frac{1}{720}h^4-\frac{1}{30240}h^6+O(h^8),\label{(3.17)}
\end{eqnarray}
where $B_n$ are Bernoulli numbers.
\end{theorem}

\emph{Proof.} We rewrite the expression (\ref{(2.7)}) as follows
\begin{eqnarray*}
\|\mathring{\ell}\|^2&=&-\sum\limits_{\beta=0}^N\mathring{C}_1[\beta]\left(\sum\limits_{\gamma=0}^N\mathring{C}_1[\gamma]G_2''(h\beta-h\gamma)-F(h\beta)\right)
+\sum\limits_{\beta=0}^N\mathring{C}_1[\beta]F(h\beta)\\
&&+\sum\limits_{\beta=0}^N\sum\limits_{\gamma=0}^NC_0[\beta]C_0[\gamma]G_2(h\beta-h\gamma)-2\sum\limits_{\beta=0}^NC_0[\beta]\int\limits_0^1G_2(x-h\beta)dx+
\int\limits_0^1\int\limits_0^1G_2(x-y)dxdy,
\end{eqnarray*}
where $F(h\beta)$ is defined by (\ref{(3.3)}).\\
Hence, taking into account (\ref{(3.16)}), we have
\begin{equation}\label{(3.18)}
\|\mathring{\ell}\|^2=A_1+A_2-2A_3+A_4,
\end{equation}
here
$$
\begin{array}{ll}
A_1=\sum\limits_{\beta=0}^N\mathring{C}_1[\beta]\ F(h\beta),& A_2=\sum\limits_{\beta=0}^N\sum\limits_{\gamma=0}^NC_0[\beta]C_0[\gamma]G_2(h\beta-h\gamma),\\
A_3=\sum\limits_{\beta=0}^NC_0[\beta]\int\limits_0^1G_2(x-h\beta)dx,& A_4=\int\limits_0^1\int\limits_0^1G_2(x-y)dxdy.
\end{array}
$$
Now we need the following sums which are obtained by using (\ref{(C0)}) and (\ref{(3.14)})
\begin{equation}\label{(3.19)}
\begin{array}{l}
\sum\limits_{\beta=0}^NC_0[\beta]=1,\ \ \sum\limits_{\beta=0}^NC_0[\beta](h\beta)=\frac{1}{2},\ \
\sum\limits_{\beta=0}^NC_0[\beta](h\beta)^2=\frac{h^2}{6}+\frac{1}{3},\\
\sum\limits_{\beta=0}^N\mathring{C}_1[\beta]e^{-h\beta}=\frac{1}{2}(1-e^{-1})\left(\frac{h(e^h+1)}{e^h-1}-2\right),\ \
\sum\limits_{\beta=0}^N\mathring{C}_1[\beta]e^{h\beta}=\frac{1}{2}(1-e)\left(\frac{h(e^h+1)}{e^h-1}-2\right).
\end{array}
\end{equation}
Taking into account (\ref{(3.9)}) and (\ref{(2.6)}), using (\ref{(3.19)}) for $A_1$, $A_2$, $A_3$ and $A_4$ we get
$$
\begin{array}{ll}
A_1=\frac{1-e^2}{8e}\left(\frac{h(e^h+1)}{e^h-1}-2\right)^2,& A_2=\frac{h^2(e^h+1)^2(e^2-1)}{8e(e^h-1)^2}-\frac{h^2+2}{12}-\frac{h(e^h+1)}{2(e^h-1)},\\
A_3=\frac{h(e^h+1)(e^2-1)}{4e(e^h-1)},& A_4=\frac{e^2-1}{2e}-\frac{7}{6}.
\end{array}
$$
Further, putting the last equalities to (\ref{(3.18)}) and after some simplifications we have
$$
\|\mathring{\ell}\|^2=1-\frac{1}{2}h+\frac{1}{12}h^2-\frac{h}{e^h-1}.
$$
Hence, using well known formula $\frac{x}{e^x-1}=\sum\limits_{n=0}^{\infty}\frac{B_n}{n!}x^n,$ $|x|<2\pi$, we get (\ref{(3.17)}).

Theorem \ref{Thm3.3} is proved\hfill $\Box$

\medskip

\textbf{Remark 2.} It should be noted that optimality of the classical Euler-Maclaurin was proved and the error of this quadrature formula was calculated in $L_2^{(m)}$, where $L_2^{(m)}$ is the space of functions which are square integrable with $m$-th generalized derivative  (see, for instance, \cite{CatCom,Schoen65,ShadHayNur13}). In particular, when $m=2$ from Corollary 5.1 of the work \cite{ShadHayNur13} we get optimality
of the Euler-Maclaurin formula
\begin{equation}\label{(3.20)}
\int_0^1\varphi(x)dx\cong h\left(\frac{1}{2}\varphi(0)+\varphi(h)+\varphi(2h)+...+\varphi(h(N-1))+\frac{1}{2}\varphi(1)\right)+\frac{h^2}{12}(\varphi'(0)-\varphi'(1))
\end{equation}
in the space $L_2^{(2)}(0,1)$.
Furthermore for the square of the norm of the error functional the following is valid
\begin{equation}\label{(3.21)}
\|\mathring{\ell}|L_2^{(2)*}(0,1)\|^2=\frac{h^4}{720}.
\end{equation}
Comparison of equalities (\ref{(3.17)}) and (\ref{(3.21)}) shows that the error of the optimal quadrature formula of the form (1.1)
on the space $W_2^{(2,1)}(0,1)$ is less than the error of the Euler-Maclaurin quadrature formula (\ref{(3.20)}) on the space $L_2^{(2)}(0,1)$.

\section*{Acknowledgments}

This work has been done while A.R.Hayotov was visiting Department of Mathematical
Sciences at KAIST, Daejeon, Republic of Korea. A.R.Hayotov is very grateful to professor Chang-Ock Lee and his research group for hospitality.
A.R. Hayotov's work was supported by the 'Korea Foundation for Advanced Studies'/'Chey Institute for Advanced Studies' International Scholar Exchange Fellowship for academic year of 2018-2019


\begin{thebibliography}{99}

\bibitem{Ahlb67} J.H. Ahlberg, E.N. Nilson, J.L. Walsh, The Theory of Splines and Their Applications, Academic Press, New
York -- London, 1967.

\bibitem{IBab}  I. Babu\v{s}ka, Optimal quadrature formulas (Russian), Dokladi
Akad. Nauk SSSR.  149 (1963) 227--229.

\bibitem{BlaCom} P. Blaga,  Gh. Coman, Some problems on optimal
quadrature, Stud. Univ. Babe\c{s}-Bolyai Math. 52, no. 4 (2007)
21--44.


\bibitem{CatCom}  T. Catina\c{s}, Gh. Coman, Optimal quadrature formulas based on the $\phi$-function method, Stud. Univ. Babe\c{s}-Bolyai Math. 51, no. 1
(2006) 49--64.


\bibitem{HayMilShad10} A.R. Hayotov, G.V. Milovanovi\'{c}, Kh.M. Shadimetov, On an optimal quadrature formula in the sense
of Sard. Numerical Algorithms, v.57, no. 4, (2011) 487-510.

\bibitem{HayMilShad15} A.R. Hayotov, G.V. Milovanovi\'{c}, Kh.M. Shadimetov,  Optimal quadratures in the sense of Sard in a Hilbert space.
Applied Mathematics and Computation, 259 (2015) 637-653.

\bibitem{Koh} P. K\"{o}hler, On the weights of Sard's quadrature
formulas, Calcolo, 25 (1988) 169--186.

\bibitem{FLan}  F. Lanzara, On optimal quadrature formulae, J. Ineq. Appl. 5 (2000)  201--225.

\bibitem{Mic74} C.A. Micchelli, Best quadrature formulas at equally spaced nodes, J. Math. Anal. Appl. 47 (1974) 232-249.

\bibitem{Nik88} S.M. Nikol'skii, Quadrature Formulas, Nauka, Moscow, 1988.(in Russian).

\bibitem{Sard}  A. Sard, Best approximate integration formulas; best approximation formulas, Amer. J. Math. 71 (1949) 80--91.

\bibitem{Schoen65} I.J. Schoenberg, On monosplines of least deviation and best quadrature formulae, J. Soc. Indust.
Appl. Math. Ser. B Numer. Anal. 2 (1965) 144-170.

\bibitem{SchSil} I.J. Schoenberg,  S.D. Silliman, On semicardinal
quadrature formulae.  Math. Comp. 28  (1974) 483--497.

\bibitem{ShadHay04} Kh.M. Shadimetov, A.R. Hayotov, Construction of the discrete analogue of the differential operator $\frac{d^{2m}}{dx^{2m}}-\frac{d^{2m-2}}{dx^{2m-2}}$, Uzbek mathematical journal, 2004, no.2, pp. 85-95.

\bibitem{ShadHay11}  Kh.M. Shadimetov, A.R. Hayotov, Optimal quadrature formulas with positive coefficients in $L_2^{(m)}(0,1)$ space,
J. Comput. Appl. Math. 235 (2011) 1114--1128.

\bibitem{ShadHay14} Kh.M. Shadimetov, A.R. Hayotov, Optimal quadrature for\-mulas in the sense of Sard in $W_2^{(m,m-1)}$ space, Calcolo
51 (2014) 211--243.

\bibitem{ShadHayNur13} Kh.M. Shadimetov, A.R. Hayotov, F.A. Nuraliev, On an optimal quadrature formula in Sobolev space $L_2^{(m)}(0,1)$,
J. Comput. Appl. Math. 243 (2013) 91--112.

\bibitem{ShadHayNur16} Kh.M. Shadimetov, A.R. Hayotov, F.A. Nuraliev, Optimal quadrature formulas of Euler-Maclaurin type,
Applied Mathematics and Computation 276 (2016) 340--355.

\bibitem{ShadNur18} Kh.M. Shadimetov, F.A. Nuraliev, Optimal formulas of numerical integration with derivatives in Sobolev space,
Journal of Siberian Federal University. Math. and Phys. 2018, 11 (6) 764-775.

\bibitem{Sobolev06}  S.L. Sobolev, The coefficients of optimal quadrature formulas, Selected Works of S.L. Sobolev,  Springer,  (2006) 561--566.

\bibitem{Sobolev74} S.L. Sobolev, Introduction to the  Theory of Cubature Formulas (Russian), Nauka, Moscow,  1974.

\bibitem{SobVas} S.L. Sobolev,  V.L. Vaskevich, The Theory of Cubature Formulas, Kluwer Academic Publishers Group, Dordrecht, 1997.

\bibitem{Zag}  F.Ya. Zagirova,  On construction of optimal quadrature
formulas with equal spaced nodes (Russian). Novosibirsk  (1982),
28 p. (Preprint No. 25, Institute of Mathematics SD of AS of USSR)

\bibitem{Zhen81} A.A. Zhensikbaev, Monosplines of minimal norm and the best quadrature formulas, Uspekhi Mat. Nauk. 1981, 36, 107--159.


\end{thebibliography}
\end{document}